\documentclass[11pt, a4paper, UKenglish]{article}

\usepackage[UKenglish]{babel}
\usepackage{amssymb, amsmath, textcomp, amsthm}
\usepackage{comment}
\usepackage{bbm,dsfont}
\usepackage{german}
\usepackage[pdftex]{graphicx}
\usepackage[utf8]{inputenc}
\usepackage{authblk}
\numberwithin{equation}{section}
\title{Hilbert transform along measurable vector fields constant on Lipschitz curves: $L^2$ boundedness}
\author{Shaoming Guo}
\date{}

\def\R{\mathbf{R}}
\def\N{\mathbf{N}}

\def\Z{\mathbf{Z}}

\def\lesim{\lesssim}

\def\begineq{\begin{equation}}
\def\endeq{\end{equation}}

\theoremstyle{plain}
\newtheorem{thm}{Theorem}[section]
\newtheorem{prop}[thm]{Proposition}
\newtheorem{lem}[thm]{Lemma}
\newtheorem{cor}[thm]{Corollary}
\newtheorem{defi}[thm]{Definition}

\newtheorem{rem}[thm]{Remark}

\newtheorem*{conj*}{Conjecture}
\newtheorem*{question*}{Question}
\newtheorem*{openproblem*}{Open Problem}

\begin{document}
\maketitle

\begin{abstract}
We prove the $L^2$ boundedness of the directional Hilbert transform 
in the plane relative to measurable vector fields which are 
constant on suitable Lipschitz curves.
\end{abstract}

\let\thefootnote\relax\footnote{Date: \date{\today}}

\section{Introduction and statement of the main result}

Consider the directional Hilbert transform in the plane defined 
for a fixed direction $v=(1,u)$ as 
\begineq
H_v f(x,y):=p.v. \int_{\R} f(x -t, y- u t)\frac{dt}t
\endeq
for any test function $f$.
By the dilation symmetry, the length of the vector $v$ is irrelevant for 
the value of $H_v$, which explains our normalization of the first component.
By an application of Fubini's theorem and the $L^p$ bounds for the classical 
Hilbert transform one obtains a priori $L^p$ bounds for $H_v$. 
On the other hand, the corresponding maximal operator $\sup_u |H_v f(x,y)|$ 
for varying directions is well known to not satisfy 
any a priori $L^p$ bounds.\\

In \cite{BT}, Bateman and Thiele proved that
\begin{equation}\label{bt}
\| \sup_{u\in \R} \| H_v f(x,y)\|_{L^p(y)} \|_{L^p(x)}\le C_p \|f\|_p
\end{equation}
for the range $3/2<p<\infty$. This result is half of a maximal bound
since the supremum is squeezed between the two $L^p$ norms on the left hand 
side, compared to being completely inside or outside as in our previous 
remarks.  The case $p=2$ of (\ref{bt}) goes back to Coifman 
and El Kohen (see \cite{BT} for the detailed discussion), who noticed that a Fourier transform in the $y$ direction makes $(\ref{bt})$ for $p=2$ equivalent to $L^2$ bounds for Carleson's 
operator. 

Estimate (\ref{bt}) highlights a bi-parameter structure of the directional
Hilbert transform. The bi-parameter structure arises since the kernel is a 
tensor product between a Hilbert kernel in direction v and a Dirac delta 
distribution in the perpendicular direction.

If one considers the linearized maximal operator 
\begin{equation}\label{linmax}
H_{v} f(x,y):=p.v. \int f(x -t, y-{u}(x,y) t)\frac{dt}t 
\end{equation}
for some function ${u}$, then inequality (\ref{bt}) can be rephrased 
as a bound for the linearized maximal operator under the assumption 
that $u$ is constant on every vertical line $x=x_0$ for all $x_0\in \R$. Such vector fields $v$ of the form $(1, u(x_0))$ for some measurable function $u:\R\to \R$ are call one-variable vector fields in \cite{BT}.\\

The purpose of the present paper is to relax 
this rigid assumption on $u$, and prove an analogue of (\ref{bt}) 
for vector fields which are constant along suitable families of 
Lipschitz curves. To formulate such a result, we perturb (\ref{bt}) by 
a bi-Lipschitz horizontal distortion, that is
\begin{equation}\label{bilip}(x,y)\to (g(x,y),y)
\end{equation}
with
\begin{equation}\label{hordis}
(x'-x)/a_0\le g(x',y)-g(x,y))\le a_0 (x'-x)
\end{equation}
for every $x<x'$ and every $y$, 
such that the transformation (\ref{bilip}) maps vertical lines to near vertical 
Lipschitz curves:
\begin{equation}
\label{liplines} 
|g(x,y)-g(x,y')|\le b_0|y'-y|
\end{equation}
for all $x, y, y'$. These two conditions can be rephrased
as 
\begineq
1/a_0\le \partial_1g\le a_0 \text{ and } |\partial_2g|\le b_0 \text{ a.e.}
\endeq

Under these assumptions, $L^p$ norms are distorted boundedly under the 
transformation (\ref{bilip}). Namely,
(\ref{hordis}) implies for every $y$ that 
\begineq\label{MM1.8}
a_0^{-1}\|f(x,y)\|_{L^p(x)}^p  \le \|f(g(x,y),y)\|_{L^p(x)}^p\le a_0 \|f(x,y)\|_{L^p(x)}^p
\endeq
and we may integrate this in $y$ direction to obtain equivalence
of $L^p$ norms in the plane.  
Hence the change of measure is not the main point of the following theorem, 
but rather the effect of the transformation on the linearizing function 
$u$, which is now constant along the family of Lipschitz curves which are
the images of the lines $x=x_0$ under the map (\ref{bilip}). 
\begin{thm}[Main Theorem]
Let $g:\R^2\to \R$ satisfy assumption (\ref{hordis}) for some $a_0$ and 
assumption (\ref{liplines}) for some $b_0$. Then for any $c_0\in (0, 1)$, we have
\begin{equation}\label{sg}
\| \sup_{|u|\le c_0/b_0} \| H_v f(g(x,y),y)\|_{L^2(y)} \|_{L^2(x)}\le C \|f\|_2.
\end{equation}
Here $C$ is a constant depending only on $a_0$ and $c_0$.
\end{thm}
\begin{rem}
The constant $C$ is independent of $b_0$ due to the unisotropic scaling symmetry $x\to x, y\to \lambda y$.
\end{rem}

In view of the implicit function theorem (see \cite{AS} for recent developments), our result covers a large class of vector fields which are of the critical Lipschitz regularity. Indeed, it implies the following

\begin{cor}\label{corollary1.2}
For a measurable unit vector field $v_0:\R^2\to S^1$, suppose that 
\vspace{2mm}

i) there exists a bi-Lipschitz map $g_0: \R^2\to \R^2$ s.t.
\begineq\label{E1.8}
v_0(g_0(x, y)) \text{ is constant in } y;
\endeq 

ii) there exists $d_0>0$ s.t. $\forall x\in \R$,
\begineq\label{scaling}
\angle (\partial_2 g_0(x, y), \pm v_0(g_0(x, y))) \ge d_0 \text{ for } y\text{-}a.e. \text{ in }\R.
\endeq 

\noindent Then the associated Hilbert transform, which is defined as
\begineq
H_{v_0}f(x, y):=\int_{\R}f((x, y)-t v_0(x, y))dt/t,
\endeq 
is bounded in $L^2$, with the operator norm depending only on $d_0$ and the bi-Lipschitz norm of $g_0$.
\end{cor}
\begin{rem}\label{NNremark1.4}
The structure theorem for Lipschitz functions by Azzam and Schul in \cite{AS} states exactly that any Lipschitz function $u:\R^2\to \R$ (any Lipschitz unit vector field $v_0$ in our case) can be precomposed with a bi-Lipschitz function $g_0:\R^2\to \R^2$ such that $u\circ g_0$ is Lipschitz in the first coordinate and constant in the second coordinate, when restricted to a ``large'' portion of the domain.
\end{rem}
\begin{rem}
Without the assumption that $d_0>0$, the operator $H_{v_0}$ might be unbounded in $L^p$ for any $p>1$. The counter-example is based on the Besicovitch-Kakeya set construction, which will be discussed at the end of the proof of the corollary.
\end{rem}

To our knowledge, this is the first result in the context of the directional
Hilbert transform with a Lipschitz assumption in the hypothesis. Lipschitz 
regularity is critical for the directional Hilbert transform as we will 
elaborate shortly.

To use the assumption that $v$ is constant along Lipschitz curves,
we apply an adapted Littlewood-Paley theory along the level lines of
$v$. This is a refinement of the analysis of Coifman and El Kohen who 
use a Fourier transform in the $y$ variable and the analysis of Bateman 
and Thiele who use a classical Littlewood-Paley theory in the $y$ variable.
This adapted Littlewood-Paley theory is the main novelty of the present paper.
It is in the spririt of prior work on the Cauchy integral on Lipschitz curves, 
for example \cite{CJS} , but it differs from this classical theme in that it 
is more of bi-parameter type as it is governed by a whole fibration into 
Lipschitz curves. 
We crucially use Jones' beta numbers as a tool to control the adapted 
Littlewood-Paley theory. To our knowledge this is also the first use of Jones' 
beta numbers in the context of the directional Hilbert transform.\\

In this paper we focus on the case $L^2$, since our goal here is to 
highlight the use of the adapted Littlewood-Paley theory and Jones' 
beta numbers in the technically most simple case. We expect to address 
the more general case $L^p$ with a range of $p$ as in the Bateman-Thiele 
theorem in forthcoming work.

While Coifman and El Kohen use the difficult bounds on Carleson's operator
as a black box, Bateman and Thiele have to unravel this black box following
the work of Lacey and Li \cite{LL}-\cite{LL1}  and use time-frequency analysis to 
prove bounds for a suitable generalization of Carleson's operator. Luckily, 
in the present work we do not have to delve into time-frequency analysis
as we can largely recycle the work of Bateman and Thiele 
for this aspect of the argument.

An upper bound such as $|u|\le c_0/b_0$ is necessary in our theorem.
By  a limiting argument we may recover the theorem of Bateman and Thiele, using 
the scaling to tighten the Lipschitz constant $b_0$ at the same time as 
relaxing the condition $|u|\le c_0/b_0$.

An interesting open question remains whether the same holds true for $c_0=1$. We do not know of a soft argument to achieve this 
relaxation, as the norm bounds arising from our argument grow unlimited as 
$c_0$ approaches $1$. This question suggests itself for further study.\\

Part of our motivation is a long history of studies of the linearized maximal 
operator (\ref{linmax}) under various assumptions on the linearizing function 
$u$. If one truncates $(\ref{linmax})$ as
$$H_{v,\epsilon_0} f(x,y):=p.v. \int_{-\epsilon_0}^{\epsilon_0} f(x -t, y-{u}(x,y) t)\frac{dt}t\ ,$$
then it is reasonable to ask for pure regularity assumptions on $u$ to obtain
boundedness of $H_{v,\epsilon_0}$. It is known that Lipschitz regularity of $u$
is critical, since a counterexample in \cite{LL1} based on a construction of the Besicovitch-Kakeya set 
shows that no bounds are possible for $C^\alpha$ regularity with $\alpha<1$.
However, it remains open whether Lipschitz regularity suffices for bounds for
$H_{v,\epsilon_0}$. On the regularity scale, the only known result is for 
real analytic vector fields $v$ by Stein and Street in \cite{SS}.
 A prior partial result in this direction appears in \cite{CNSW}.

It is our hope that our result corners some of the difficulties of 
approaching Lipschitz regularity in the classical problem. Further 
substantial progress (including the case $c_0=1$) is likely to use Lipschitz regularity not only of 
the level curves of $u$ but also of $u$ itself across the level curves. For example, one possibility would be to cut the plane into different pieces by the theorem of Azzam and Schul stated in Remark \ref{NNremark1.4}, and to analyze each piece separately by the Main Theorem in the present paper. We leave this for the future study.\\

Related to the directional Hilbert transform and thus additional motivation
for the present work is the directional maximal operator
\begin{equation}
M_{v,\epsilon_0} f(x):= \sup_{0<\epsilon<\epsilon_0}\frac{1}{2\epsilon} \int_{-\epsilon}^{\epsilon} |f(x-t, y- u(x,y) t )|dt
\end{equation}
which arises for example in Lebesgue type differentiation questions and 
has an even longer history of interest than the directional Hilbert transform.  
Hilbert transforms and maximal operators share many features, in particular 
they have the same scaling and thus share the same potential $L^p$ bounds.
The maximal operator is in some ways easier as it is positive and does not
have a singular kernel. For example, bounds for the maximal operator under the
assumption of real analytic vector fields were proved much earlier by Bourgain in \cite{Bo2}.

An instance of bounds satisfied by the maximal operator but not the Hilbert transform arises when one restricts the range
of the function $u$ instead of the regularity.  For certain sets of directions characterized by Bateman in \cite{Ba5} there are bounds for the maximal operator (for example for the set of lacunary directions), while
Karagulyan proves in \cite{Kara} that no such bounds are possible for the Hilbert transform.

On the other hand, the Hilbert transform is easier in some other aspects, 
most notably it is a linear operator. For example, bounds for the bilinear
Hilbert transform mapping into $L^1$ were known (Lacey and Thiele \cite{LT1}, \cite{LT2})
before the corresponding maximal operator bounds (Lacey \cite{La}), due to the fact that orthogonality between different tiles is preserved under the Hilbert transform but not the maximal operator. In particular we 
do not know at the moment whether the analogue of our main theorem holds for 
the directional maximal operator. This may be an interesting subject for 
further investigation.\\

{\bf Outline of paper:} in Section 2 we will prove Corollary \ref{corollary1.2} by reducing it to the Main Theorem. The reduction will also be used later in the proof of the Main Theorem.

In Section 3 we will state the strategy of the proof for the Main Theorem. As it appears that our result is a Lipschitz perturbation of the one by Bateman and Thiele, this turns out also to be the case for the proof: if we denote by $P_k$ a Littlewood-Paley operator in the $y$-variable, the main observation in Bateman and Thiele's proof is that $H_v$ commutes with $P_k$. In our case, this is no longer true. However, we can make use of an adapted version of the Littlewood-Paley projection operator $\tilde{P}_k$ (see Definition \ref{definition1}) to partially recover the orthogonality. We split the operator $H_v$ into a main term and a commutator term
\begineq
\sum_{k\in \Z} H_v P_k(f)=\sum_{k\in \Z}(H_v P_k(f)-\tilde{P}_k H_v P_k(f)+\tilde{P}_k H_v P_k(f)).
\endeq

The boundedness of the main term $\sum_{k\in \Z}\tilde{P}_k H_v P_k(f)$ is essentially due to Lacey and Li \cite{LL1}, with conditionality on certain maximal operator estimate. In Section 4 we modify Bateman's argument in \cite{Ba} and \cite{Ba2} to the case of vector fields constant on Lipschitz curves and remove the conditionality on that maximal operator. 

The main novelty is the boundedness of the commutator term 
\begineq
\sum_{k\in \Z}(H_v P_k(f)-\tilde{P}_k H_v P_k(f)),
\endeq 
which will be presented in Section 5. To achieve this, we will view Lipschitz curves as perturbations of straight lines and use Jones' beta number condition for Lipschitz curves and the Carleson embedding theorem to control the commutator. Here we shall emphasis again that the commutator estimate is free of time-frequency analysis.\\

{\bf Notations:} throughout this paper, we will write $x\ll y$ to mean that $x\le y/10$, $x\lesim y$ to mean that there exists a universal constant $C$ s.t. $x\le C y$, and $x\sim y$ to mean that $x\lesim y$ and $y\lesim x$.  $\mathbbm{1}_E$ will always denote the characteristic function of the set $E$.\\

{\bf Acknowledgements.} This work is done under the supervision of Prof. Christoph Thiele, to whom the author would like to express his most sincere gratitude, for suggesting him such an interesting topic, for sharing with him his deep insight into this problem, and for numerous suggestions on the exposition of this paper. 



\section{Proof of Corollary \ref{corollary1.2}}

In this section we prove Corollary \ref{corollary1.2}, by reducing it to the Main Theorem. Some parts of the reduction will also be used in the proof of the Main Theorem in the rest of the paper.\\

We first divide the unit circle $S^1$ into $N$ arcs of equal length, with the angle of each arc being $2\pi/N$. Choose 
\begineq
N> 6\pi/d_0
\endeq
s.t. $2\pi/N< d_0/3$. Denote these arcs as $\Omega_1, \Omega_2, ..., \Omega_N$. For each $\Omega_i$, define
\begin{displaymath}
H_{v_0, \Omega_i}f(x,y):=\left\{
\begin{array}{ll}
H_{v_0} f(x,y) & \textrm{if $v_0(x,y)\in \Omega_i$}\\
0		& \textrm{else}
\end{array}\right.
\end{displaymath}
If we were able to prove that $\|H_{v_0, \Omega_i}\|_{2\to 2}$ is bounded by a constant $C$ which is independent of $i\in \{1, 2, ..., N\}$, then we conclude that 
\begineq
\|H_{v_0}\|_{2\to 2}\le C N(d_0).
\endeq

Now fix one $\Omega_i$, we want to show the boundedness of $H_{v_0, \Omega_i}$. Choose a new coordinate such that the $x$-axis passes through $\Omega_i$ and bisects it. Then all the vectors in $\Omega_i$ form an angle less than $d_0/6$ with the $x$-axis. As we assume that \begineq
\angle(\partial_2 g_0, \pm v_0(g_0))\ge d_0>0,
\endeq 
we see that the vector $\partial_2 g_0$ forms an angle less than $\frac{\pi-d_0}{2}$ with the $y$-axis.

Renormalize the unit vector $v_0$ such that the first component is 1, i.e. write $v_0=(1, u_0)$, then by the fact that $v_0$ forms an angle less than $d_0/6$ with the $x$-axis, we obtain 
\begineq
|u_0|\le \tan(d_0/6).
\endeq

Next we construct the Lipschitz function $g$ in the Main Theorem from the bi-Lipschitz map $g_0$, and the coordinate we will use here is still the one associated to $\Omega_i$ as above. Under this linear change of variables, we know that $g_0$ is still bi-Lipschitz. We renormalize the bi-Lipschitz map in such a way that 
\begineq
g_0(x,0)=(x,0), \forall x\in \R.
\endeq
Fix $x\in \R$, the map $g_0$, when restricted on the vertical line $\{(x,y): y\in \R\}$, is still bi-Lipschitz. We denote by $\Gamma_x$ the image of this bi-Lipschitz map, i.e.
\begineq
\Gamma_x:=\{g_0(x, y): y\in \R\}.
\endeq 
Define the function $g$ by the following relation
\begineq
(g(x,y),y)=g_0(x, y\prime),
\endeq
for some $y\prime$. By the fact that $g_0$ is bi-Lipschitz, we know that such $y\prime$ exists and is unique.\\

From the above construction and the fact that $\partial_2 g_0$ forms an angle less than $\frac{\pi-d_0}{2}$ with the $y$-axis, we see easily that
\begineq\label{EE2.9}
|g(x,y_1)-g(x,y_2)|\le \cot(d_0/2) |y_1-y_2|, \forall x, y_1, y_2\in \R.
\endeq

Hence what is left is to show that condition \eqref{hordis} is also satisfied with a constant $a_0$ depending only on $d_0$ and the bi-Lipschitz constant of $g_0$. One side of the equivalence $(x_1-x_2)/a_0\le g(x_1, y)-g(x_2, y)$ is quite clear from the picture below: the bi-Lipschitz map $g_0$ sends the points $P, Q$ to $(g(x_1, y), y), (g(x_2, y), y)$ separately, then by definition of bi-Lipschitz map, there exists constant $a_0$ s.t. 
\begineq
g(x_1, y)-g(x_2, y)\ge \frac{1}{a_0}|P-Q|\ge \frac{1}{a_0}(x_1-x_2).
\endeq
\begin{minipage}[b]{10.2cm}
\includegraphics[width=10cm]{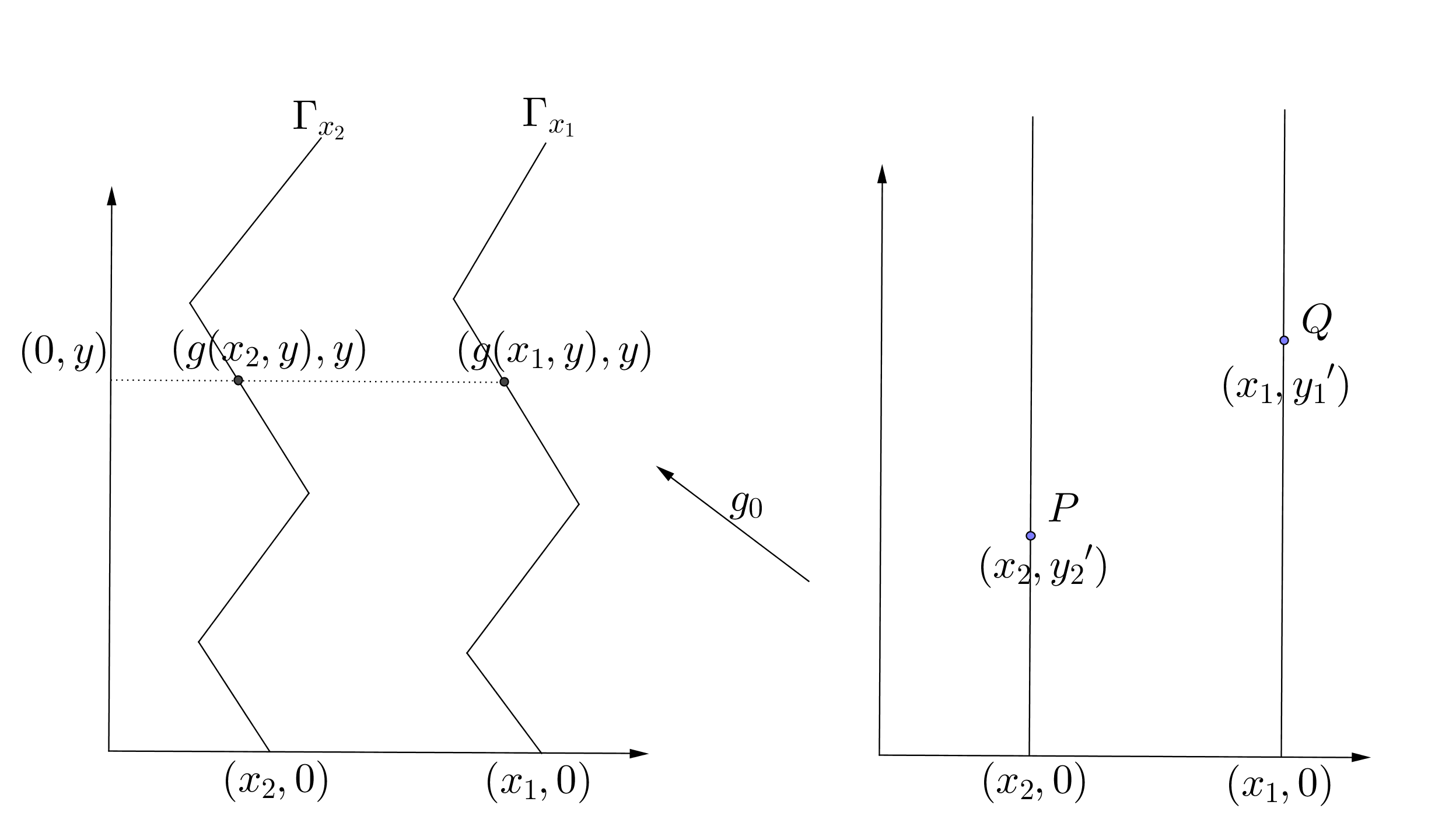}
\end{minipage}

\noindent For the other side, we argue by contradiction. If for any $M\in \N$ large, there exists $x_1, x_2, y\in \R$ s.t. 
\begineq
g(x_1, y)-g(x_2, y)\ge M(x_1-x_2),
\endeq
then together with \eqref{EE2.9}, this implies that 
\begineq
dist(K, \Gamma_{x_1})\ge M\sin(d_0/2)(x_1-x_2).
\endeq
But this is not allowed as by the definition of the bi-Lipschitz map $g_0$ and the Lipschitz function $g$, $dist(K, \Gamma_{x_1})$ must be comparable to $|x_1-x_2|$.

So far, we have verified all the conditions in the Main Theorem with
\begineq
b_0=\cot(d_0/2) \text{ and } c_0=\tan(d_0/6)/\cot(d_0/2)<1.
\endeq
Hence we can apply the Main Theorem to obtain the boundedness of $H_{v_0,\Omega_i}$.\\

In the end, as claimed in the corollary, we still need to show the blowing up of the operator norm without the assumption that $d_0>0$. Indeed, the blowing up happens in $L^p$ not only for $p=2$ but also for all $p>1$. For the range $p\le 2$, the counter example is simply a Knapp example: take the function $f(x)=\mathbbm{1}_{B_1(0)}(x)$, let $\Gamma$ be the upper cone which forms an angle less than $\frac{\pi}{4}$ with the vertical axis. First define the vector field $v(x)=\frac{x}{|x|}$ for $x\in \Gamma\setminus B_1(0)$, then extend the definition to the whole plane properly such that $v$ satisfies the condition \eqref{E1.8}. It is then easy to see that
\begineq
|H_v f(x)|\sim \frac{1}{|x|}, \forall x\in \Gamma\setminus B_1(0),
\endeq
which does not belong to $L^p(\R^2)$ for $p\le 2.$ For the range $p>2$, the counter example is given by the standard Besicovitch-Kakeya set construction, which can be found in \cite{Ba2} and \cite{LL1}.

\section{Strategy of the proof of the Main Theorem}

If we linearize the maximal operator in the Main Theorem, what we need to prove turns to be the following
\begineq\label{MM3.1}
\|\int_{\R}f(g(x,y)-t, y-t u(x))dt/t\|_2\lesim \|f\|_2,
\endeq
where $u:\R\to \R$ is a measurable function such that $\|u\|_{\infty}\le c_0/b_0$. The change of coordinates 
\begineq
(x, y)\to (g(x, y), y)
\endeq
in \eqref{bilip} also changes the measure on the plane. However, we still want to use the original Lebesgue measure for the Littlewood-Paley decomposition. Hence we invert \eqref{bilip} and denote the inversion by
\begineq\label{MM3.3}
(x, y)\to (P(x, y), y),
\endeq
where ``P'' stands for ``projection'', and the reason of calling so can be illustrated by the following picture:

\begin{minipage}[b]{15.2cm}
\includegraphics[width=10cm]{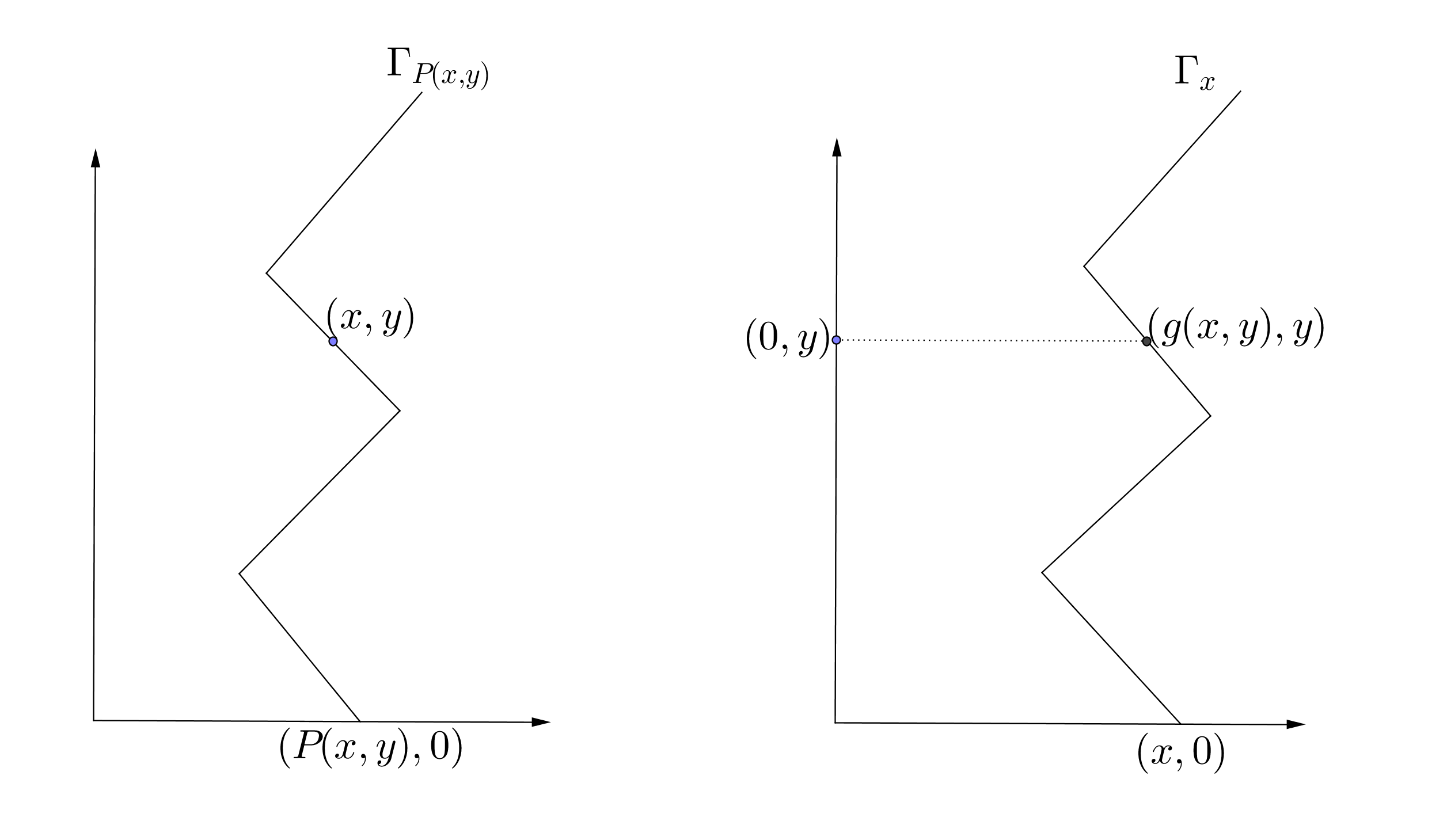}
\end{minipage}

The change of coordinates in \eqref{MM3.3} turns the estimate \eqref{MM3.1} into the following equivalent form 
\begineq
\|\int_{\R}f(x-t, y-t u(P(x,y)))dt/t\|_2 \lesim \|f\|_2.
\endeq
Moreover, we will denote
\begineq\label{MM3.5}
H_v f(x, y):=\int_{\R}f(x-t, y-t u(P(x,y)))dt/t.
\endeq
In the rest of the paper, we want to make the convention that whenever $H_v$ appears, it denotes the Hilbert transform along the vector field $v(x, y)=(1, u(P(x, y)))$, which is the above \eqref{MM3.5}, just to distinct it from the various $H_v$ that have appeared in the introduction.\\

To prove the above estimate, we first make several reductions: by the unisotropic scaling 
\begineq
x\to x, y\to \lambda y,
\endeq
we can w.l.o.g. assume that $b_0=10^{-2}$. By a similar cutting and pasting argument to that in the proof of Corollary \ref{corollary1.2}, we can assume that $c_0\ll 10^{-2}$, i.e. the vector field $v$ is of the form $(1, u)$ with $|u|\ll 1$. \\

Now we start the detail of the proof. It was already observed in Bateman $\cite{Ba2}$ that under the assumption $|u|\ll 1$, we can w.l.o.g. assume that supp $\hat{f}$ lies in a two-ended cone which forms an angle less than $\pi/4$ with the vertical axis, as for functions $f$ with frequency supported on the cone near the horizontal axis, we have that
\begineq
H_vf(x, y)=H_{(1, 0)}f(x, y),
\endeq 
which is the Hilbert transform along the constant vector field $(1, 0)$. But $H_{(1,0)}$ is bounded by Fubini's theorem and the $L^2$ boundedness of the Hilbert transform.

For the frequencies outside the cone near the horizontal axis, the proof consists of two steps. In the first step we will prove the boundedness of $H_v$ when acting on functions with frequency supported in one single annulus. To be precise, let $\Gamma$ be the cone which forms an angle less than $\pi/4$ with the vertical axis, $\Pi_{\Gamma}$ be the projection operator on $\Gamma$, i.e.
\begineq
\Pi_{\Gamma} f:=\mathcal{F}^{-1} \mathbbm{1}_{\Gamma} \mathcal{F}f,
\endeq
where $\mathcal{F}$ stands for the Fourier transform and $\mathcal{F}^{-1}$ the inverse transform. Let $P_k$ be the $k$-th Littlewood-Paley projection operator in the vertical direction, as we are always concerned with the frequency in $\Gamma$, later for simplicity we will just write $P_k$ instead of $P_k \Pi_{\Gamma}$ for short. Then what we will prove first is 
\begin{prop}\label{thm1}
Under the same assumptions as in the Main Theorem, we have for $p\in (1, \infty)$ that
\begineq
\|H_v P_k (f)\|_p\lesim \|P_k(f)\|_p,
\endeq
with the constant being independent of $k\in \Z$.
\end{prop}

In order to prove the boundedness of $H_v$, we need to put all the frequency pieces together. In the case of $C^{1+\alpha}$ vector fields for any $\alpha>0$, Lacey and Li's idea in \cite{LL1} is to prove the almost orthogonality between different frequency annuli. In the case where the vector field is constant along vertical lines, an important observation in the paper of Bateman and Thiele is that $H_v$ and $P_k$ commute, which then makes it possible to apply a Littlewood-Paley square function estimate. 

In our case Bateman and Thiele's observation is no longer true. We need to take into account that the vector field is constant along Lipschitz curves, which gives rise to an adapted Littlewood-Paley projection operator (the following Definition \ref{definition1}). 


Before defining this operator, we first need to make some preparation. Fix one $\tilde{x}\in \R$, take the curve $\Gamma_{\tilde{x}}$ which passes through $(\tilde{x}, 0)$, recall that $\Gamma_{\tilde{x}}$ is given by the set $\{(g(\tilde{x}, \tilde{y}),\tilde{y}): \tilde{y}\in \R\}$, where $g$ is the Lipschitz function in the Main Theorem. By the definition of the operator $H_v$ we know that the vector field $v$ is equal to the constant vector $(1, u(\tilde{x}))$ along $\Gamma_{\tilde{x}}$. Change the coordinate s.t. the horizontal $x^{\prime}$-axis is parallel to $(1, u(\tilde{x}))$. The following lemma says that in the new coordinate, the curve $\Gamma_{\tilde{x}}$ can still be realized as the graph of a Lipschitz function.

\begin{lem}\label{lem3.2}
For any fixed $\tilde{x}\in \R$, there exists a Lipschitz function $x^{\prime}=g_{\tilde{x}}(y^{\prime})$ s.t. $\Gamma_{\tilde{x}}$ can be reparametrized as $\{(g_{\tilde{x}}(y^{\prime}), y^{\prime}): y^{\prime}\in \R\}$. Moreover, we have that $\|g_{\tilde{x}}\|_{Lip}\le \frac{1+b_0}{1-b_0}$, where $b_0$ is the constant in the Main Theorem.
\end{lem}
\noindent {\bf Proof of Lemma \ref{lem3.2}:} denote by $\theta$ the angle between the vector $(1, u(\tilde{x}))$ and the $x$-axis as in the picture below.\\
\begin{minipage}[b]{10.2cm}
\includegraphics[width=10cm]{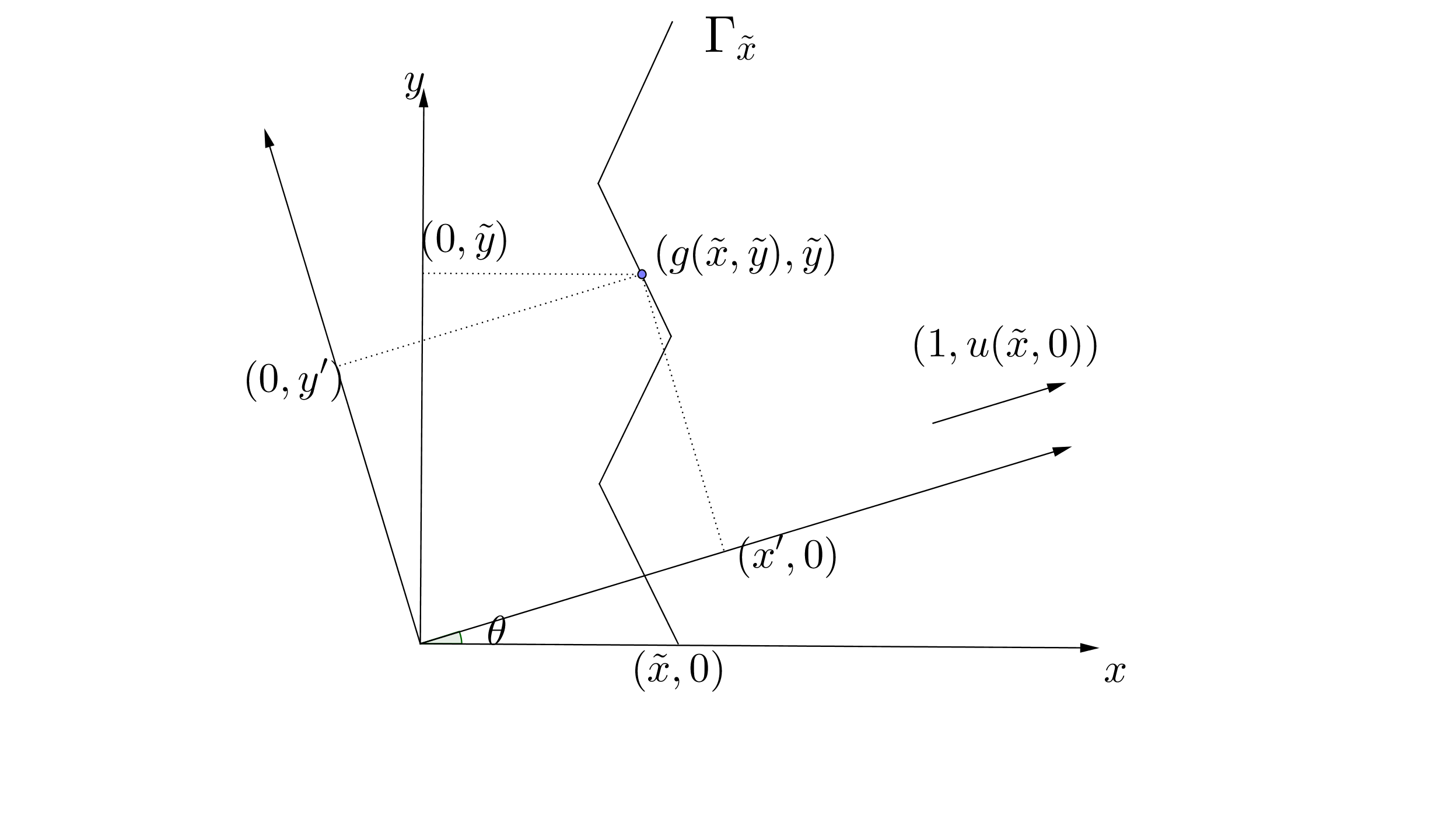}
\end{minipage}

\noindent The new coordinate of the point $(g(\tilde{x}, \tilde{y}), \tilde{y})$ will be given by
\begineq\label{E3.4}
(x^{\prime}, y^{\prime})=(\tilde{y} \sin \theta+g(\tilde{x}, \tilde{y})\frac{1+\sin^2\theta}{\cos \theta}, \tilde{y} \cos \theta-g(\tilde{x}, \tilde{y})\sin\theta).
\endeq
Look at the identity for the second component
\begineq
y^{\prime}=\tilde{y} \cos \theta-g(\tilde{x}, \tilde{y})\sin\theta,
\endeq
we want to solve $\tilde{y}$ by $y^{\prime}$ by using the implicit function theorem. As
\begineq
\frac{dy^{\prime}}{d\tilde{y}}=\cos \theta-\frac{\partial g}{\partial \tilde{y}} \sin \theta,
\endeq
by the fact that $|u|\ll 1$ and $|\frac{\partial g}{\partial \tilde{y}}|\le b_0\le 10^{-2}$, we obtain that 
\begineq
\frac{1-b_0}{\sqrt{2}}\le \frac{dy^{\prime}}{d\tilde{y}}\le \frac{1+b_0}{\sqrt{2}},
\endeq 
from which it is clear that the implicit function theorem is applicable.

After solving $\tilde{y}$ by $y^{\prime}$, we just need to substitute $\tilde{y}$ into the identity for the first component in \eqref{E3.4}, which is 
\begineq
x^{\prime}=\tilde{y} \sin \theta+g(\tilde{x}, \tilde{y})\frac{1+\sin^2\theta}{\cos \theta},
\endeq
to get an implicit expression of $x^{\prime}$ in terms of $y^{\prime}$, which we will denote as $x^{\prime}=g_{\tilde{x}}(y^{\prime})$.

To estimate the Lipschitz norm of the function $g_{\tilde{x}}$, we just need to observe that when doing the above change of variables, we have rotated the axis by an angle $\theta$ which satisfies $|\theta|\le \pi/4$. Together with the fact that $|\frac{\partial g}{\partial \tilde{y}}|\le b_0$, we can then derive that 
\begineq
|\frac{\partial g_{\tilde{x}}}{\partial y^{\prime}}|\le \frac{1+b_0}{1-b_0},
\endeq 
which finishes the proof of Lemma \ref{lem3.2}. $\Box$

\begin{defi}(adapted Littlewood-Paley projection)\label{definition1} Select a Schwartz function $\psi_0$ with support on $[\frac{1}{2},\frac{5}{2}]\cup [-\frac{5}{2},-\frac{1}{2}]$ such that 
\begineq
\sum_{k\in \Z} \psi_0(2^{-k}t)=1, \forall t\neq 0.
\endeq
For $f:\R^2\to \R$, for every fixed $\tilde{x}\in \R$, define the adapted (one dimensional) Littlewood-Paley projection on $\Gamma_{\tilde{x}}$ by
\begineq
\tilde{P}_k(f)(x^{\prime}, y^{\prime}):=\int_{\R}f(g_{\tilde{x}}(z), z) \check{\psi}_k(y^{\prime}-z)dz=P_k(\tilde{f})(y^{\prime}),
\endeq
where $(x^{\prime}, y^{\prime})=(g_{\tilde{x}}(y^{\prime}), y^{\prime})$ denotes one point in $\Gamma_{\tilde{x}}$, $\psi_k(\cdot):=\psi_0(2^{-k} \cdot)$ and we use $\tilde{f}(\cdot)$ to denote the function $f(g_{\tilde{x}}(\cdot), \cdot)$, and $P_k$ the one dimensional Littlewood-Paley projection operator.
%
%
\end{defi}

Now it is instructive to regard the Lipschitz curves as perturbation of the straight lines, or equivalently, to think that $H_v P_k f$ still has frequency supported near the $k$-th frequency band, which has already been used by Lacey and Li in their almost orthogonality estimate for $C^{1+\alpha}$ vector fields in \cite{LL1}. We then subtract the term $\tilde{P}_k H_v P_k(f)$ from $H_v P_k(f)$, and estimate the commutator. 

To be precise, we first write
\begineq\label{E3.5}
\sum_{k} H_v P_k(f)=\sum_{k}(H_v P_k(f)-\tilde{P}_k H_v P_k(f)+\tilde{P}_k H_v P_k(f)),
\endeq
then by the triangle inequality, we have
\begineq
\|\sum_{k} H_v P_k(f)\|_2\lesim\|\sum_{k}(H_v P_k(f)-\tilde{P}_k H_v P_k(f))\|_2+\|\sum_{k}\tilde{P}_k H_v P_k(f)\|_2.
\endeq
We call the second term the main term, and the first term the commutator term. The $L^2$ boundedness of the main term will follow from orthogonality argument, which is the following adapted Littlewood-Paley theorem.

\begin{lem}\label{lemma3}
For $p\in (1,+\infty)$, we have the following variants of the Littlewood-Paley estimates:
\begineq\label{EE3.9}
\| (\sum_{k\in \Z}|\tilde{P}_k(f)|^2)^{1/2}\|_p\sim \|f\|_p,
\endeq
\begineq\label{E3.10}
\| (\sum_{k\in \Z}|\tilde{P}_k^{*}(f)|^2)^{1/2}\|_p\sim \|f\|_p,
\endeq
with constants depending only on $a_0$.
\end{lem}
\noindent {\bf Proof of Lemma \ref{lemma3}:} in \eqref{MM1.8} from the introduction, we have already explained the following coarea formula:
\begineq\label{E3.11}
\int_{\R^2}|f(x, y)|dxdy\sim\int_{\R}[\int_{\Gamma_{\tilde{x}}} |f|ds_{\tilde{x}}]d\tilde{x}.
\endeq

\noindent We apply this formula to the left hand side of \eqref{EE3.9} to obtain
\begineq\label{E3.9}
\| (\sum_{k\in \Z}|\tilde{P}_k(f)|^2)^{1/2}\|_p^p \sim \int_{\R} \int_{\Gamma_{\tilde{x}}} (\sum_{k\in \Z} |\tilde{P}_k (f)|^2)^{p/2}ds_{\tilde{x}} d\tilde{x}.
\endeq
For every fixed $\tilde{x}$, by Definition \ref{definition1}, the right hand side of \eqref{E3.9} turns to 
\begineq
\int_{\R}[\int_{\R} (\sum_{k} |P_k(\tilde{f}_{\tilde{x}})(y^{\prime})|^2)^{p/2}dy^{\prime}]d\tilde{x},
\endeq
where $\tilde{f}_{\tilde{x}}(y^{\prime})=f(g_{\tilde{x}}(y^{\prime}), y^{\prime})$. Then the classical Littlewood-Paley theory applies and we can bound the last expression by
\begineq
\int_{\R} \|f\|_{L^p(\Gamma_{\tilde{x}})}^p d\tilde{x} \lesim \|f\|_{L^p}^p.
\endeq

For the boundedness of the adjoint operator, it suffices to prove that 
\begineq\label{E3.12}
\sum_{k\in \Z} \langle \tilde{P}_k^*(f), f_k\rangle \lesim \|f\|_{L^p} \|(\sum_{k\in \Z} |f_k|^2)^{1/2}\|_{L^{p^{\prime}}}.
\endeq
First by linearity and H\"older's inequality, we derive 
\begineq
\sum_{k\in \Z} \langle \tilde{P}_k^*(f), f_k\rangle=\langle f, \sum_{k\in \Z} \tilde{P}_k (f_k)\rangle \lesim \|f\|_{L^p} \|\sum_{k\in \Z} \tilde{P}_{k} (f_k)\|_{L^{p^{\prime}}}.
\endeq
Applying the coarea formula \eqref{E3.11}, we obtain
\begineq
\|\sum_{k\in \Z} \tilde{P}_{k} (f_k)\|_{L^{p^{\prime}}}\sim (\int_{\R}(\int_{\Gamma_{\tilde{x}}}|\sum_{k\in \Z}\tilde{P}_k(f_k)|^{p'}ds_{\tilde{x}})d\tilde{x})^{1/p'}.
\endeq
By the Definition \ref{definition1}, for every fixed $\tilde{x}\in \R$, the inner integration in the last expression turns to 
\begineq\label{E3.25}
\int_{\R}|\sum_{k\in \Z}P_k(\tilde{f}_{k, \tilde{x}})(y^{\prime})|^{p'}dy^{\prime},
\endeq
where $\tilde{f}_{k, \tilde{x}}(y^{\prime}):=f_k(g_{\tilde{x}}(y^{\prime}), y^{\prime})$. Now the classical Littlewood-Paley theory applies and we bound the term in \eqref{E3.25} by 
\begineq\label{E3.26}
\int_{\R}(\sum_{k\in \Z}|\tilde{f}_{k, \tilde{x}}(y^{\prime})|^2)^{p'/2}dy^{\prime}\lesim \int_{\Gamma_{\tilde{x}}}(\sum_{k\in \Z}|f_k|^2)^{p'/2}ds_{\tilde{x}}\lesim \|(\sum_{k\in \Z}|f_k|^2)^{1/2}\|_{L^{p'}(\Gamma_{\tilde{x}})}^{p'}.
\endeq
Then to prove \eqref{E3.12}, we just need to integrate $d\tilde{x}$ in \eqref{E3.26} and apply the coarea formula \eqref{E3.11} to derive 
\begin{align*}
\|\sum_{k\in \Z} \tilde{P}_{k} (f_k)\|_{L^{p^{\prime}}}& \lesim (\int_{\R}\|(\sum_{k\in \Z}|f_k|^2)^{1/2}\|_{L^{p'}(\Gamma_{\tilde{x}})}^{p'}d\tilde{x})^{1/p'}\\
		&\lesim \|(\sum_{k\in \Z} |f_k|^2)^{1/2}\|_{L^{p^{\prime}}}.
\end{align*}
Thus we have finished the proof of Lemma \ref{lemma3}.$\Box$\\

\noindent Now we will show how to prove the $L^2$ boundedness of the main term by Lemma \ref{lemma3} and Proposition \ref{thm1}: first by duality, we have
\begin{align*}
\|\sum_{k} \tilde{P}_k H_v P_k(f)\|_2 &= \sup_{\|g\|_2=1}|\langle \sum_{k}\tilde{P}_k H_v P_k(f), g\rangle|\\
		&= \sup_{\|g\|_2=1}|\langle\sum_{k}H_v P_k(f), \tilde{P}_k^{*}(g)\rangle|.
\end{align*}
Applying the Cauchy-Schwartz inequality and H\"older's inequality, we can bound the last term by
\begineq
\sup_{\|g\|_2=1} \|(\sum_{k}|H_v P_k(f)|^2)^{1/2}\|_2 \|(\sum_{k} |\tilde{P}_k^{*}(g)|^2)^{1/2}\|_{2}.
\endeq
For the former term, Proposition \ref{thm1} implies that
\begin{align*}
 \|(\sum_{k}|H_v P_k(f)|^2)^{1/2}\|_2& \le (\sum_{k\in Z}\|H_v P_k(f)\|_2^2)^{1/2}\\
	&\lesim (\sum_{k\in \Z}\|P_k(f)\|_2^2)^{1/2}\lesim \|f\|_2.
\end{align*}
For the latter term, Lemma \ref{lemma3} implies that 
\begineq
 \|(\sum_{k} |\tilde{P}_k^{*}(g)|^2)^{1/2}\|_{2} \lesim \|g\|_2.
\endeq
Thus we have proved the $L^2$ boundedness the main term, modulo Proposition \ref{thm1}.\\  

As the second step, we will prove the $L^2$ boundedness of the commutator, which is
\begineq
\|\sum_{k}(H_v P_k(f)-\tilde{P}_k H_v P_k(f))\|_2 \lesim \|f\|_2.
\endeq
To do this, we first split the operator $H_v$ into a dyadic sum: select a Schwartz function $\psi_0$ such that $\psi_0$ is supported on $[\frac{1}{2},\frac{5}{2}]$, let 
\begineq
\psi_l(t):=\psi_0(2^{-l} t),
\endeq
by choosing $\psi_0$ properly, we can construct a partition of unity for $\R^+$, i.e. 
\begineq
\mathbbm{1}_{(0,\infty)}=\sum_{l\in \Z} \psi_l.
\endeq
Let 
\begineq\label{E4.24}
H_{l}h(x,y):= \int \check{\psi}_l(t)h(x-t, y-tu(P(x,y)))dt,
\endeq
then the operator $H_v$ can be decomposed into the sum
\begineq
H_v=-\mathds{1}+ 2 \sum_{l\in \Z} H_{l}.
\endeq
Hence to bound the commutator, it is equivalent to bound the following
\begineq
\sum_{k\in \Z} \sum_{l\in \Z}(H_{l} P_k f- \tilde{P}_k H_{l} P_k f).
\endeq
Notice that by definition, $H_{l} P_k f$ vanishes for $l>k$, which simplifies the last expression to 
\begineq
\sum_{l\ge 0} \sum_{k\in \Z} (H_{k-l} P_k f- \tilde{P}_k H_{k-l} P_k f).
\endeq
By the triangle inequality, it suffices to prove
\begin{prop}\label{thm2}
Under the same assumption as in the Main Theorem, there exists $\gamma>0$ such that
\begineq
\|\sum_{k\in \Z} (H_{k-l} P_k f- \tilde{P}_k H_{k-l} P_k f)\|_2\lesim 2^{-\gamma l}\|f\|_2,
\endeq
with the constant independent of $l\in \N.$
\end{prop}
So far, we have reduced the proof of the Main Theorem to that of Proposition \ref{thm1} and Proposition \ref{thm2}, which we will present separately in the following sections.

\section{Boundedness of the Lipschitz-Kakeya maximal function and proof of Proposition \ref{thm1}}

Lacey and Li in their prominent work \cite{LL1} have reduced the $L^2$ boundedness of the operator $H_{v, \epsilon_0}$ to the boundedness of one new operator they introduced, which is the so called Lipschitz-Kakeya maximal operator. As soon as this operator is bounded, we can then repeat the argument in \cite{LL1} to obtain Proposition \ref{thm1} as a corollary. 

Here we follow \cite{Ba2}, where a slightly different version of the Lipschitz-Kakeya maximal operator is used, see the following Lemma \ref{lem4.8}.  The only place in \cite{Ba2} where the one-variable vector field plays a special role is Lemma 6.2 in page 1037. Hence to prove Proposition \ref{thm1}, we just need to replace this lemma by Lemma \ref{lem4.8}, and leave the rest of the argument unchanged.\\

In this section we make an observation that both the boundedness of the Lipschitz-Kakeya maximal operator (Corollary \ref{lemma1}) and its variant (Lemma \ref{lem4.8}) can be proved by adapting Bateman's argument in \cite{Ba} to our case where the vector fields are constant only on Lipschitz curves.

Before defining the Lipschitz-Kakeya maximal operator, we first need to introduce several definitions.

\begin{defi}(popularity)
For a rectangle $R\subset \R^2$, with $l(R)$ its length, $w(R)$ its width, we define its uncertainty interval $EX(R)\subset \R$ to be the interval of width $w(R)/l(R)$ and centered at slope($R$). Then the popularity of the rectangle $R$ is defined to be 
\begineq
pop_R:=|\{(x,y)\in \R^2: u(P(x,y))\in EX(R)\}|/|R|.
\endeq
\end{defi}
\begin{defi}
Given two rectangles $R_1$ and $R_2$ in $\R^2$, we write $R_1\le R_2$ whenever $R_1\subset C R_2$ and $EX(R_2)\subset EX(R_1)$, where $C$ is some properly chosed large constant, and $C R_2$ is the rectangle with the same center as $R_2$ but dilated by the factor $C$. 
\end{defi}

Denote $\mathcal{R}_{\delta, \omega}:=\{R\in \mathcal{R}: slope(R)\in [-1,1], pop_R\ge \delta, w(R)=\omega\}$, where $\mathcal{R}$ is the collection of all the rectangles in $\R^2$. Then the Lipschitz-Kakeya maximal function is defined as
\begineq\label{maxoperator}
M_{\mathcal{R}_{\delta, \omega}}(f)(x):=\sup_{x\in R\in \mathcal{R}_{\delta, \omega}}\frac{1}{|R|}\int_R |f|
\endeq

\begin{lem}\label{lem4.8}
Let $u$ and $P$ be the functions given in the definition of the operator $H_v$ in \eqref{MM3.5}. Suppose $\mathcal{R}_0$ is a collection of pairwise incomparable (under ``$ \le $'') rectangles of uniform width such that for each $R\in \mathcal{R}_0$, we have
\begineq
\frac{|(u\circ P)^{-1}(EX(R))\cap R|}{R}\ge \delta, (\text{i.e. } pop_{R}\ge \delta)
\endeq
and 
\begineq
\frac{1}{|R|}\int_R \mathbbm{1}_F \ge \lambda.
\endeq
Then for each $p>1$, 
\begineq
\sum_{R\in \mathcal{R}_0}|R| \lesim \frac{|F|}{\delta \lambda^p}.
\endeq
\end{lem}

The same covering lemma argument as in Lemma 3.1 \cite{Ba} shows the boundedness of Lacey and Li's Lipschitz-Kakeya maximal operator as a corollary of Lemma \ref{lem4.8}.
\begin{cor}\label{lemma1}
For all $p\in (1,\infty)$ we have the following bound
\begineq
\|M_{\mathcal{R}_{\delta, \omega}}\|_{L^p\to L^P}\le C(p,a_0)\frac{1}{\delta}
\endeq
\end{cor}

\noindent {\bf Proof of Lemma \ref{lem4.8}:} the proof is essentially due to Bateman \cite{Ba}. Most of the argument in \cite{Ba} remains, with just one minor modification in order to adapt to the family of Lipschitz curves on with the vector field is constant.

\begin{defi}\label{definition4}(rectangles adapted to the vector field)
For a rectangle $R\in \mathcal{R}_{\delta, \omega}$, with its two long sides lying on the parallel lines $y=kx+b_1$ and $y=kx+b_2$ for some $k\in [-1, 1]$ and $b_1, b_2\in \R$, define $\tilde{R}$ to be the adapted version of $R$, which is given by the set
\begineq
\{(x, y): P(x, y)\in P(R)\}\bigcap \{(x, kx+b):x\in \R, b\in [b_1, b_2]\},
\endeq
where $P$ is the projection operator in \eqref{MM3.3}.
\end{defi}

What we need to do is just to replace the rectangles $R$ in \cite{Ba} by $\tilde{R}$, and observe that the two key quantities--length and popularity of rectangles-- are both preserved under the projection operator $P$ up to a constant depending on the constant $a_0$ in the Main Theorem. Hence we leave out the details and refer to \cite{Ba}.$\Box$

\section{Proof of Proposition \ref{thm2}}

This section consists of three subsections. In the first subsection we will introduce some notations, most of which we adopt from Bateman's paper \cite{Ba2}, with minor changes for our purpose. In the second we will use Jones' beta numbers and the Carleson embedding theorem to prove Proposition \ref{thm2}, modulo one crucial lemma which will be presented afterwards in the third subsection.

\subsection{Discretization}
The content of this subsection is basically taken from Bateman's paper \cite{Ba2}, with minor changes as we are now dealing with all frequencies instead of one single frequency annulus.\\

{\bf Discretizing the functions:} Fix $l\ge 0$, we write $\mathcal{D}_l$ as the collection of the dyadic intervals of length $2^{-l}$ contained in $[-2,2]$. Fix a smooth positive function $\beta:\R\to \R$ s.t. 
\begineq
\beta(x)=1, \forall |x|\le 1; \beta(x)=0, \forall |x|\ge 2.
\endeq
Also choose $\beta$ such that $\sqrt{\beta}$ is a smooth function. Then fix an integer $c$(whose exact value is unimportant), for each $\omega\in \mathcal{D}_l$, define
\begineq
\beta_{\omega}(x)= \beta(2^{l+c}(x-c_{\omega_1})),
\endeq
where $\omega_1$ is the right half of $\omega$ and $c_{\omega_1}$ is its center.

Define
\begineq
\beta_l(x)=\sum_{\omega\in \mathcal{D}_l} \beta_{\omega}(x),
\endeq
note that 
\begineq
\beta_l(x+2^{-l})=\beta_l(x), \forall x\in [-2, 2-2^{-l}].
\endeq
Define 
\begineq
\gamma_l=\frac{1}{2}\int_{-1}^1 \beta_l(x+t)dt,
\endeq
because of the above periodicity, we know that $\gamma_l$ is constant for $x\in [-1,1]$, independent of $l$. Say $\gamma_l(x)=\delta>0$, hence 
\begineq
\frac{1}{\delta}\gamma_l(x) \mathds{1}_{[-1,1]}(x)=\mathds{1}_{[-1,1]}(x).
\endeq
Define another multiplier $\tilde{\beta}:\R\to \R$ with support in $[\frac{1}{2}, \frac{5}{2}]$ and $\tilde{\beta}(x)=1$ for $x\in [1,2]$. We define the corresponding multiplier on $\R^2$:
\begin{align*}
& \hat{m}_{k,\omega}(\xi, \eta)=\tilde{\beta}(2^{-k}\eta)\beta_{\omega}(\frac{\xi}{\eta})\\
& \hat{m}_{k,l,t}(\xi, \eta)=\tilde{\beta}(2^{-k}\eta)\beta_l(t+\frac{\xi}{\eta})\\
& \hat{m}_{k,l}(\xi, \eta)=\tilde{\beta}(2^{-k}\eta)\gamma_l(\frac{\xi}{\eta})
\end{align*}
Then what we need to bound can be written as 
\begin{align*}
\|\sum_{k\in \Z}\sum_{l\in \Z}H_l P_k(f)\|_p &=\|\int_{-1}^1 \sum_{k\in \Z}\sum_{l\ge 0} H_{k-l}(\frac{1}{\delta}m_{k,l}*f) dt\|_p\\
			&\le \int_{-1}^1 \|\sum_{k\in \Z} \sum_{l\ge 0} H_{k-l} (\frac{1}{\delta}m_{k,l,t}*f)\|_p dt,
\end{align*}
where the terms $H_l P_k$ for $l>k$ in the sum vanish as explained before.

So it suffices to prove a uniform bound on $t\in [-1, 1]$, w.l.o.g. we will just consider the case $t=0$, which is 
\begineq
\sum_{k\in \Z} \sum_{l\ge 0} H_{k-l}(m_{k,l,0}*f)=\sum_{k\in \Z} \sum_{l\ge 0} H_{k-l}([\tilde{\beta}(2^{-k}\eta)\beta_l(\frac{\xi}{\eta})]*f).
\endeq

{\bf Constructing the tiles:} For each $k\in \Z$ and $\omega\in \mathcal{D}_l$ with $l\ge 0$, let $\mathcal{U}_{k,\omega}$ be a partition of $\R^2$ by rectangles of width $2^{-k}$ and length $2^{-k+l}$, whose long side has slope $\theta$, where $\tan \theta=-c(\omega)$, which is the center of the interval $\omega$. If $s\in \mathcal{U}_{k,\omega}$, we will write $\omega_s:=\omega$, and $\omega_{s,1}$ to be the right half of $\omega$, $\omega_{s, 2}$ the left half.


An element of $\mathcal{U}_{k, \omega}$ for some $\omega\in \mathcal{D}_l$ is called a ``tile''. Define $\varphi_{k,\omega}$ such that
\begineq
|\hat{\varphi}_{k,\omega}|^2=\hat{m}_{k,\omega},
\endeq
then $\varphi_{k,\omega}$ is smooth by our assumption on $\beta$ mentioned above.

For a tile $s\in \mathcal{U}_{k,\omega}$, define 
\begineq\label{wavelet}
\varphi_s(p):= \sqrt{|s|} \varphi_{k,\omega} (p-c(s)),
\endeq
where $c(s)$ is the center of $s$. Notice that 
\begineq
\|\varphi_s\|_2^2= \int_{\R^2} |s| \varphi_{k,\omega}^2=|s|\int_{\R^2} \hat{m}_{k,\omega}=1,
\endeq
i.e. $\varphi_s$ is $L^2$ normalized.

The construction of the tiles above by uncertainty principle is to localize the function further in space, for this purpose we need
\begin{lem}(\cite{Ba2})
\begineq
f*m_{k,\omega}(x)=\lim_{N\to \infty} \frac{1}{4 N^2} \int_{[-N,N]^2} \sum_{s\in \mathcal{U}_{k,\omega}} \langle f,\varphi_s(p+\cdot)\rangle \varphi_s(p+x)dp
\endeq
\end{lem}
\noindent The above lemma allows us to pass to the model sum
\begin{align*}
\sum_{k\in \Z} \sum_{l\ge 0} H_{k-l}(f* m_{k,l,0})= \sum_{k\in \Z}\sum_{l\ge 0}\sum_{\omega\in \mathcal{D}_l}\sum_{s\in \mathcal{U}_{k,\omega}}\langle f, \varphi_s\rangle H_{k-l}(\varphi_s),
\end{align*}
define 
\begineq
\psi_s=\psi_{-\log(length(s))},
\endeq
and 
\begineq\label{EE5.13}
\phi_s(x,y):= \int \check{\psi}_s(t)\varphi_s(x-t,y-tu(P(x,y)))dt,
\endeq
then the model sum turns to 
\begineq\label{E5.14}
\sum_{k\in \Z}\sum_{l\ge 0}\sum_{\omega\in \mathcal{D}_l}\sum_{s\in \mathcal{U}_{k,\omega}}\langle f, \varphi_s\rangle \phi_s
\endeq

\begin{lem}\label{lemma5}
we have that $\phi_s(x,y)=0$ unless $-u(P(x,y))\in \omega_{s,2}$.
\end{lem}
The proof of Lemma \ref{lemma5} is by the Plancherel theorem, we just need to observe that the frequency support of $\psi_s$ and $\hat{\varphi}_s$ will be disjoint at the point $(x, y)$ unless $-u(P(x,y))\in \omega_{s,2}$.

\subsection{Boundedness of the commutator and proof of Proposition \ref{thm2}}

This subsection is devoted to the proof of Proposition \ref{thm2}, which is motivated in large by the proof of the $T(b)$ theorem and the boundedness of the paraproduct, see \cite{AHMTT} and \cite{CJS} for example. 


In our case, unlike Bateman and Thiele's proof for the one-variable vector fields, it's no longer true that $H_v P_k f$ still has frequency in the $k$-th annulus. In order to get enough orthogonality for the term $H_v P_k f$ to apply the Littlewood-Paley theory, we need to subtract the term $H_v P_k f- \tilde{P}_k H_v P_k f$, which should be viewed as a family of paraproducts.\\

We proceed with the details of the proof. If we expand the summation on the left hand side of Proposition \ref{thm2} with \eqref{E5.14}, what we need to bound can be rewritten as 
\begineq
 \|\sum_k  \sum_{\omega\in \mathcal{D}_l}\sum_{s\in \mathcal{U}_{k,\omega}} \langle f, \varphi_s\rangle (\phi_s-\tilde{P}_k \phi_s)\|_2\lesim 2^{-\gamma l}\|f\|_2.
\endeq

In order to use the orthogonality of different wave packets, we will prove the $L^2$ bound  for the dual operator, which is
\begineq
\sum_k  \sum_{\omega\in \mathcal{D}_l}\sum_{s\in \mathcal{U}_{k,\omega}} \langle h, \phi_s-\tilde{P}_k \phi_s\rangle \varphi_s.
\endeq
Notice that for $s_1\in \mathcal{U}_{k_1, \omega_1}$ and $s_2\in \mathcal{U}_{k_2, \omega_2}$ with $(k_1, \omega_1)\neq(k_1, \omega_2)$,  we have
\begineq\label{orthogonality}
\langle \varphi_{s_1}, \varphi_{s_2}\rangle=0
\endeq
by the definition of the wavelet function $\varphi_s$ in \eqref{wavelet}. Also if we know that $s_1, s_2$ are in the same $\mathcal{U}_{k,\omega}$, for some $k$ and $\omega$, then we can find $m,n\in \Z$ s.t.
\begineq
c(s_2)=c(s_1)+(m\cdot l(s_1), n\cdot w(s_1))
\endeq
where $c(s)$ is the center of the tile $s$, $l(s)$ its length and $w(s)$ its width. Then by the non-stationary phase method we know for any $N\in \N$, there exists a constant $C_N$ depending only on $N$ s.t.
\begineq\label{almostorthogonality}
|\langle \varphi_{s_1}, \varphi_{s_2}\rangle|\le \frac{C_N}{<|m|+|n|>^N}.
\endeq 
Here we want to make a remark that the exact value of $N$ is not important, it just denotes some large number which might vary from line to line if we use the same notation later.\\

\noindent Applying the above two estimates \eqref{orthogonality} \eqref{almostorthogonality}, we obtain
\begin{align*}
& \|\sum_{k}\sum_{\omega\in \mathcal{D}_{l}}\sum_{s\in \mathcal{U}_{k,\omega}}\langle h, \phi_s-\tilde{P}_k \phi_s\rangle\varphi_s\|_{2}^{2}\\
& = \sum_{k}\sum_{\omega\in \mathcal{D}_{l}}\sum_{s_{1}\in \mathcal{U}_{k,\omega}}\sum_{s_{2}\in \mathcal{U}_{k,\omega}}\langle h, \phi_{s_{1}}-\tilde{P}_k \phi_{s_{1}}\rangle \langle\varphi_{s_{1}},\varphi_{s_{2}}\rangle \langle h, \phi_{s_{2}}-\tilde{P}_k \phi_{s_{2}}\rangle.
\end{align*}
As we know for any $s_{1},s_{2}\in \mathcal{U}_{k, \omega}$ there exists $m,n\in \Z$ s.t.
\begineq\label{E1.66}
c(s_2)=c(s_1)+(m\cdot l(s_1), n\cdot w(s_1)),
\endeq
the above sum can be rewritten as
\begineq
\sum_{m,n\in \Z} \sum_{k\in \Z}\sum_{\omega\in \mathcal{D}_{l}}\sum_{s_{1}\in \mathcal{U}_{k,\omega}}\langle h, \phi_{s_{1}}-\tilde{P}_k \phi_{s_{1}}\rangle \langle\varphi_{s_{1}},\varphi_{s_{2}}\rangle \langle h, \phi_{s_{2}}-\tilde{P}_k \phi_{s_{2}}\rangle
\endeq
with $s_{1},s_{2}$ satisfying the relation \eqref{E1.66}. \\

\noindent Now fix $m,n\in \Z$, by the estimate in \eqref{almostorthogonality}, we know that 
\begin{align*}
& \sum_{k}\sum_{\omega\in \mathcal{D}_{l}}\sum_{s_{1}\in \mathcal{U}_{k,\omega}}|\langle h, \phi_{s_{1}}-\tilde{P}_k \phi_{s_{1}}\rangle \langle\varphi_{s_{1}},\varphi_{s_{2}}\rangle \langle h, \phi_{s_{2}}-\tilde{P}_k \phi_{s_{2}}\rangle|\\
& \lesim \frac{1}{<|m|+|n|>^{N}}\sum_{k}\sum_{\omega\in \mathcal{D}_{l}}\sum_{s_{1}\in \mathcal{U}_{k,\omega}}|\langle h, \phi_{s_{1}}-\tilde{P}_k \phi_{s_{1}}\rangle \langle h, \phi_{s_{2}}-\tilde{P}_k \phi_{s_{2}}\rangle|,
\end{align*}
by the Cauchy-Schwartz inequality, the last term is bounded by 
\begineq
\frac{1}{<|m|+|n|>^{N}}\sum_{k}\sum_{\omega\in \mathcal{D}_{l}}\sum_{s\in \mathcal{U}_{k,\omega}}|\langle h, \phi_{s}-\tilde{P}_k \phi_{s}\rangle|^{2},
\endeq
then it suffices to prove that
\begineq\label{E5.23}
\sum_k \sum_{\omega\in \mathcal{D}_{l}}\sum_{s\in \mathcal{U}_{k,\omega}}\langle h, \phi_s-\tilde{P}_k \phi_s\rangle^2 \lesim 2^{-\gamma l}\|h\|_2^2.
\endeq

%
%

First to estimate every single term $\langle h, \phi_s-\tilde{P}_k \phi_s\rangle$ for a fixed tile $s$:  denote $s_{m,n}$ to be the shift of $s$ by $(m,n)$ units, i.e. 
\begineq
s_{m,n}:=\{(x,y)\in \R^2: (x-m\cdot l(s), y-n\cdot w(s))\in s\},
\endeq
then by the triangle inequality we know that
\begineq
|\langle h, \phi_s-\tilde{P}_k \phi_s\rangle|\le  \sum_{m,n\in \Z}|\int_{s_{m,n}} h\cdot (\phi_s-\tilde{P}_k \phi_s)dydx|.
\endeq
Recall that in Definition \ref{definition4} we use $\tilde{R}$ to denote the adapted version of the rectangle $R$ to the family of Lipschitz curves, then clearly $\tilde{s}_{m,n}\supset s_{m,n}$. Thus
\begineq
|\langle h, \phi_s-\tilde{P}_k \phi_s\rangle|\le  \sum_{m,n\in \Z}|\int_{\tilde{s}_{m,n}} h\cdot (\phi_s-\tilde{P}_k \phi_s)dydx|.
\endeq
By the coarea formula \eqref{E3.11}, we obtain
\begin{align*}
|\langle h, \phi_s-\tilde{P}_k \phi_s\rangle|& \le  \sum_{m,n\in \Z}|\int_{\tilde{s}_{m,n}} h\cdot (\phi_s-\tilde{P}_k \phi_s)dydx|\\
		&\lesssim \sum_{m,n\in \Z} \int_{P(s_{m,n})}\int_{\Gamma_x\cap \tilde{s}_{m,n}}|h\cdot (\phi_s-\tilde{P}_k \phi_s)|ds_x dx,
\end{align*}
where $ds_x$ stands for the arc length measure of the Lipschitz curve $\Gamma_x$. 

Now for the inner integration along the curve $\Gamma_x$, we do the same change of coordinates and the same parametrization of $\Gamma_x$ as in Definition \ref{definition1}, i.e. we choose the coordinates s.t. the horizontal axis is parallel to $(1, u(x))$, and represent the curve $\Gamma_x$ by the Lipschitz function $g_x(\cdot)$. If we let $J(x, s_{m,n})$ denote the projection of $\Gamma_x\cap \tilde{s}_{m,n}$ on the new vertical axis, the last expression becomes
\begineq\label{EE5.27}
\sum_{m,n\in \Z} \int_{P(s_{m,n})}\int_{J(x, s_{m,n})}|h(g_x(y), y) (\phi_s(g_x(y), y) - P_k [\phi_s(g_x(y), y)])|dy dx.
\endeq

To bound the above term, Jones' beta number will play a crucial role.
\begin{defi}(\cite{Jo})
For a Lipschitz function $A:\R\to \R$, we first take the Calder\'on decomposition of $a(x)=A^{\prime}(x)$, which yields the representation
\begineq
a(x)=\sum_{I\text{ } dyadic} a_I \psi_I(x),
\endeq
where $\psi_I$ is some mean zero function supported on $3I$, $|\psi_I^{\prime}(x)|\le |I|^{-1}.$ For each dyadic interval $I$, let
\begineq
\alpha_I=\sum_{|J|\ge |I|} a_I \psi_J(c_I),
\endeq
where $c_I$ stands for the center of $I$, denote the ``average slope'' of the Lipschitz curve near $I$, and define the beta number
\begin{equation}
\beta_0(I):= \sup_{x\in 3 I}\frac{|A(x)-A(c_I)-\alpha_I(x-c_I)|}{|I|},
\end{equation}
and the $j_0$-th beta number
\begin{equation}
\beta_{j_0}(I):=\sup_{x\in 3 j_0 I}\frac{|A(x)-A(c_I)-\alpha_I(x-c_I)|}{|I|}.
\end{equation}
\end{defi}

For beta numbers, we have the following Carleson condition.
\begin{lem}(\cite{Jo})
For any Lipschitz function $A$, we have
\begineq
\sup_J \frac{1}{|J|} \sum_{I\subset J}\beta_0^2(I) |I| \lesim \|A\|_{Lip}^2,
\endeq
and also for any $j_0\in \N$
\begineq
\sup_J \frac{1}{|J|} \sum_{I\subset J}\beta_{j_0}^2(I) |I| \lesim j_0^3 \|A\|_{Lip}^2.
\endeq
\end{lem}

After introducing Jones' beta number, we are ready to state
\begin{lem}\label{claim1}
for $x\in P(s_{m,n})$, we have the following estimate:
\begin{align*}
& \int_{J(x, s_{m,n})}|h(g_x(y), y) (\phi_s(g_x(y), y) - P_k [\phi_s(g_x(y), y)])|dy\\
& \lesim  \sum_{j_0\in \N}\frac{2^{-3l/2}}{<|j_0|+|m|+|n|>^N}  \beta_{j_0}(x, s_{m,n}) [h]_{x, s_{m,n}} \mathbbm{1}_{\{-u(x)\in \omega_{s,2}\}}(x)
\end{align*}
where $ \beta_{j_0}(x, s_{m,n})$ is the $j_0$-th beta number for the Lipschitz curve $g_x(\cdot)$ on the interval $J(x, s_{m,n})$, $[h]_{x, s_{m,n}}$ is the average of the function $h$ on the interval $J(x, s_{m,n})$, i.e.
\begineq
[h]_{x, s_{m,n}}:=\frac{1}{w(s)}\int_{J(x, s_{m,n})}|h(g_x(y), y)|dy.
\endeq
\end{lem}

The proof of Lemma \ref{claim1} will be postponed to the next subsection. Substitute the estimate in Lemma \ref{claim1} into the estimate for the term $\langle h, \phi_s-\tilde{P}_k \phi_s\rangle$, we then have that
\begin{align*}
& |\langle h, \phi_s-\tilde{P}_k\phi_s\rangle|\\
& \lesim \sum_{m,n}\int_{P(s_{m,n})}\int_{J(x, s_{m,n})}|h(g_x(y), y) (\phi_s(g_x(y), y) - P_k [\phi_s(g_x(y), y)])|dy dx\\
&\lesim \sum_{m,n}\int_{P(s_{m,n})}\sum_{j_0\in \N}\frac{2^{-3l/2}}{<|j_0|+|m|+|n|>^N}  \beta_{j_0}(x, s_{m,n}) [h]_{x, s_{m,n}}\mathbbm{1}_{\{-u(x)\in \omega_{s,2}\}}(x)dx
\end{align*}
hence 
\begin{align*}
& \sum_k \sum_{\omega\in \mathcal{D}_{l}}\sum_{s\in \mathcal{U}_{k,\omega}}|\langle h, \phi_s-\tilde{P}_k\phi_s\rangle|^2\\
& \lesim \sum_k \sum_{\omega\in \mathcal{D}_{l}}\sum_{s\in \mathcal{U}_{k,\omega}}\sum_{m,n,j_0}\frac{2^{-3l}}{<|j_0|+|m|+|n|>^N} ...\\
& \text{ }\text{ }\text{ }\text{ }\text{ }\text{ }\text{ }\text{ }\text{ }\text{ }\text{ }\text{ }\text{ }\text{ }\text{ }\text{ }\text{ }\text{ }\text{ }\text{ }\text{ }\text{ }...|\int_{P(s_{m,n})} \beta_{j_0}(x, s_{m,n}) [h]_{x, s_{m,n}}\mathbbm{1}_{\{-u(x)\in \omega_{s,2}\}}(x)dx|^2\\
& \lesim \sum_{m,n,j_0}\frac{2^{-2l}}{<|j_0|+|m|+|n|>^N} ...\\
& \text{ }\text{ }\text{ }\text{ }\text{ }\text{ }...\sum_k \sum_{\omega\in \mathcal{D}_{l}}\sum_{s\in \mathcal{U}_{k,\omega}} w(s) \int_{P(s_{m,n})} \beta_{j_0}^2(x, s_{m,n}) [h]^2_{x, s_{m,n}}\mathbbm{1}_{\{-u(x)\in \omega_{s,2}\}}(x)dx
\end{align*}

\begin{lem}\label{lemma7}
for any fixed $x$, fixed $m,n,j_{0}$,
\begineq
\sum_k \sum_{\omega\in \mathcal{D}_{l}}\sum_{s\in \mathcal{U}_{k,\omega}} w(s) \mathbbm{1}_{P(s_{m,n})}(x)\beta_{j_0}^2(x, s_{m,n}) [h]^2_{x, s_{m,n}} \mathbbm{1}_{\{-u(x)\in \omega_{s,2}\}}(x) \lesim j_{0}^{3} \|h\|_{L^2(\Gamma_x)}^{2}
\endeq
\end{lem}
\noindent {\bf Proof of Lemma \ref{lemma7}:} this lemma is akin to the Carleson embedding theorem, as we have the following Carleson type condition
\begineq
\sup_{s_{m,n}}\frac{1}{|J(x,s_{m,n})|}\sum_{{s^\prime}_{m,n}:J(x,s^{\prime}_{m,n})\subset J(x,s_{m,n})} \beta_{j_{0}}^{2}(J(x,s^{\prime}_{m,n})) w(s^{\prime}_{m,n})\lesim j_{0}^{3} Lip^{2}(\Gamma_{x}),
\endeq
where the term $\mathbbm{1}_{\{-u(x)\in \omega_{s,2}\}}$ plays such a role that, originally there are $2^l$ groups of dyadic rectangles 
\begineq
\bigcup_k \bigcup_{\omega\in \mathcal{D}_l}\bigcup_{s\in \mathcal{U}_{k, \omega}}\{s_{m,n}\}
\endeq 
in the summation  $\sum_{k}\sum_{\omega\in \mathcal{D}_l}\sum_{s\in \mathcal{U}_{k, \omega}}$, which means that there are also $2^l$ groups of dyadic intervals 
\begineq
\bigcup_k \bigcup_{\omega\in \mathcal{D}_l}\bigcup_{s\in \mathcal{U}_{k, \omega}} \{J(x, s_{m,n})\}
\endeq 
which are the projections of the intersection of the dyadic rectangles with $\Gamma_x$ on the vertical axis, the term $\mathbbm{1}_{\{-u(x)\in \omega_{s,2}\}}$ just guarantees that there is just one such collection which has contribution, i.e. which has the right orientation in the sense of Lemma \ref{lemma5}.

Then the desired estimate will just follow from the Carleson embedding theorem, which we refer to Lemma 5.1 in \cite{AHMTT}.$\Box$\\

Continue the calculation before the above lemma:
\begin{align*}
& \sum_k \sum_{\omega\in \mathcal{D}_{l}}\sum_{s\in \mathcal{U}_{k,\omega}}|\langle h, \phi_s-\tilde{P}_k\phi_s\rangle|^2 \\
& \lesim \sum_{m,n, j_0}\frac{2^{-2l} j_0^3}{<|j_0|+|m|+|n|>^N} \int_{\R} \|h\|_{L^2(\Gamma_x)}^2 dx \lesim 2^{-2l} \|h\|_2^2.
\end{align*}

This finishes the proof for \eqref{E5.23} and then Proposition \ref{thm2} modulo Lemma \ref{claim1}, which we will present in the following subsection.\\

\subsection{Proof of Lemma \ref{claim1}}
We assume that $-u(x)\in \omega_{s, 2}$, which means the vector $(1, u(x))$ is roughly parallel to the long side of $s_{m,n}$, otherwise the left hand side in Lemma \ref{claim1} will also vanish due to Lemma \ref{lemma5}. After the change of variables in \eqref{EE5.27}, the vector $(1, u(x))$ turns to $(1, 0)$.\\

{\bf Proof by ignoring the tails:} in order to explain how Jones' $\beta$-number appears, we first sketch the proof by ignoring the tails of the wavelet functions and the tail of the kernel of the Littlewood-Paley projection operator $P_k$. 

By the above simplification, we only need to consider the case $m=n=0$. What we need to ``prove'' becomes 
\begineq\label{AA5.37}
\int_{J(x, s)} |h(g_x(y), y) (\phi_s(g_x(y), y) - P_k [\phi_s(g_x(y), y)])|dy \lesim 2^{-3l/2} \beta_0(J(x, s)) [h]_{x, s}. 
\endeq
For fixed $x$, we denote by $\tau_{x, s}y+b$ the line of ``average slope'' we picked in the definition of the beta number for the Lipschitz curve $g_x(\cdot)$ on the interval $J(x, s)$, for the sake of simplicity we assume $b=0$. Moreover, as both $x$ and $s$ are fixed, we will also just write $\tau$ instead of $\tau_{x, s}$. Then we make the crucial observation that
\begineq\label{AA5.38}
P_k[\phi_s^x(\tau y,y)]=\phi_s^x(\tau y,y),
\endeq
where 
\begineq
\phi_s^x(\tau y,y):=\int_{\R}\check{\psi}_s(t)\varphi_s(\tau y-t, y)dt,
\endeq
due to the fact that for any function $\varphi_s$ with frequency supported on the $k$-th annulus, if we restrict the function to a straight line, it will still have frequency supported on the $k$-th annulus (with one dimension less).

In comparison with the definition of $\phi_s$ in \eqref{EE5.13}, $\phi_s^x(\tau y,y)$ is defined as the Hilbert transform along the vector $(1, u(x))$ (which is (1, 0) after the change of the variables we made in Lemma \ref{lem3.2} and in the expression \eqref{EE5.27}) instead of the direction of the vector field $v$ at the point $(\tau y, y)$. 

Hence from the identity in \eqref{AA5.38} we obtain 
\begineq\label{AA5.40}
\begin{split}
& \phi_s(g_x(y),y)-P_k[\phi_s(g_x(y),y)]\\
&=\phi_s(g_x(y),y)-P_k[\phi_s(g_x(y),y)-\phi_s^x(\tau y, y)+\phi_s^x(\tau y, y)]\\
& =\phi_s(g_x(y),y)-\phi_s^x(\tau y,y)-P_k[\phi_s(g_x(y),y)-\phi_s^x(\tau y, y)].
\end{split}
\endeq
As we have also ignored the tails of the kernel of $P_k$, it is easy to see that the former and the latter terms in the last expression can essentially be handled in the same way. Hence in the following we will only consider the former term, which corresponds to the term
\begineq\label{AA5.41}
\begin{split}
& \int_{J(x, s)}|h(g_x(y), y) (\phi_s(g_x(y),y)-\phi_s^x(\tau y,y))|dy.
\end{split}
\endeq
By the definitions of $\phi_s$ and $\phi_s^x$, we have
\begineq\label{AA5.42}
\begin{split}
& |\phi_{s}(g_x(y),y)-\phi_{s}^x(\tau y, y)|\\
& = |\int_{\R} \check{\psi}_{k-l}(t)\varphi_{s}(g_x(y)-t, y)dt-\int_{\R} \check{\psi}_{k-l}(t)\varphi_{s}(\tau y-t, y)dt|\\
& = 2^{k-l} |\int_{\R}\check{\psi}_0(2^{k-l}t)\varphi_{s}(g_x(y)-t, z)dt-\int_{\R} \check{\psi}_0(2^{k-l}t)\varphi_{s}(\tau y-t, y)dt|\\
& = 2^{k-l} |\int_{\R}[\check{\psi}_0(2^{k-l}(t+g_x(y)-\tau y))-\check{\psi}_0(2^{k-l}t)]\varphi_{s}(\tau y-t, z)dt|.
\end{split}
\endeq
By the definition of the beta numbers, we have that
\begineq
|g_x(y)-\tau y|\lesim \beta_{0}(x, s) 2^{-k},
\endeq
which implies that
\begineq
|\check{\psi}_0(2^{k-l}(t+g_x(y)-\tau y))-\check{\psi}_0(2^{k-l}t)|\lesim 2^{-l} \beta_{0}(x, s)
\endeq
by the fundamental theorem. In the end, by substituting the above estimate into \eqref{AA5.42} and \eqref{AA5.41} we obtain the desired estimate \eqref{AA5.37}.\\

{\bf The full proof:} the main idea is still the same, and the difference is that we need to be more careful with the tails of the wavelet functions and the kernel of $P_k$.  

For fixed $x$, fixed $m$ and $n$, denote $\tau(x, s_{m,n})y+b$ as the line of ``average slope'' for the Lipschitz curve $g_x(\cdot)$ on the interval $J(x, s_{m,n})$, for the sake of simplicity we assume $b=0$. Then the crucial observation \eqref{AA5.38} becomes
\begineq\label{EE5.37}
P_k[\phi_s^x(\tau(x, s_{m,n})y,y)]=\phi_s^x(\tau(x, s_{m,n})y,y).
\endeq
Hence similar to \eqref{AA5.40}, we obtain from \eqref{EE5.37} that
\begin{align*}
& \phi_s(g_x(y),y)-P_k[\phi_s(g_x(y),y)]\\
& =\phi_s(g_x(y),y)-\phi_s^x(\tau(x, s_{m,n})y,y)-P_k[\phi_s(g_x(y),y)-\phi_s^x(\tau(x, s_{m,n})y, y)].
\end{align*}
Denote
\begineq
I_{s_{m,n}}=|\int_{J(x, s_{m,n})}h(g_x(y),y) \cdot (\phi_s(g_x(y),y)-\phi_s^x(\tau(x, s_{m,n})y,y))dy|
\endeq
and also 
\begineq
II_{s_{m,n}}= |\int_{J(x, s_{m,n})}h(g_x(y),y) \cdot P_k[\phi_s(g_x(y),y)-\phi_s^x(\tau(x, s_{m,n})y,y)]dy|.
\endeq

\begin{lem}\label{lemma9}
Under the above notations, for $z\in J(x, s_{m,n})+ j_0 2^{-k}$ with $j_0\in \Z$, we have the pointwise estimate
\begineq
|\phi_s(g_x(z), z)-\phi_s^x(\tau(x, s_{m,n})z, z)|\lesim \frac{\beta_{|j_0|}(x, s_{m,n}) 2^{k} 2^{-3l/2}}{<\min\{|m|+|n|, |m|+|n|-|j_0|\}>^N}.
\endeq 
\end{lem}

Let us first continue the proof of Lemma \ref{claim1}: for the first term $I_{s_{m,n}}$, we take $j_0$ in Lemma \ref{lemma9} to be zero, then
\begineq
|\phi_s(g_x(z), z)-\phi_s^x(\tau(x, s_{m,n})z, z)|\lesim \frac{\beta_{0}(x, s_{m,n}) 2^{k} 2^{-3l/2}}{<|m|+|n|>^N},
\endeq
which implies that
\begineq
I_{s_{m,n}}\lesim \frac{2^{-3l/2}}{<|m|+|n|>^N} \beta_{0}(x, s_{m,n}) [h]_{x, s_{m,n}}.
\endeq

For the second term $II_{s_{m,n}}$, by the definition of $P_k$, 
\begin{align*}
& |P_k[\phi_s(g_x(y),y)-\phi_s^x(\tau(x, s_{m,n})y,y)]|\\
&  =|\int_{\R} (\phi_s(g_x(z),z)-\phi_s^x(\tau(x, s_{m,n})z ,z))2^k \check{\psi}_0(2^k(y-z))dz|\\
& \le |\sum_{j_0\in \Z} \int_{J(x, s_{m,n})+j_0 2^{-k}} (\phi_s(g_x(z),z)-\phi_s^x(\tau(x, s_{m,n})z ,z))2^k \check{\psi}_0(2^k(y-z))dz|.
\end{align*}
For $y\in J(x, s_{m,n})$ and $z\in J(x, s_{m,n})+j_0 2^{-k}$, by the non-stationary phase method, we have that 
\begineq
|\check{\psi}_0(2^k(y-z))|\lesim \frac{1}{<j_0>^N},
\endeq
together with the estimate in Lemma \ref{lemma9}, we arrive at
\begin{align*}
& |P_k[\phi_s(g_x(y),y)-\phi_s^x(\tau(x, s_{m,n})y,y)]|\\
& \lesim \sum_{j_0\in \Z}  \frac{\beta_{|j_0|}(x, s_{m,n}) 2^{k} 2^{-3l/2}}{<\min\{|m|+|n|, |m|+|n|-|j_0|\}>^N} \frac{1}{<j_0>^N}\\
& \lesim \sum_{j_0\in \Z}  \frac{\beta_{|j_0|}(x, s_{m,n}) 2^{k} 2^{-3l/2}}{<|m|+|n|+|j_0|>^N}.
\end{align*}
Substitute the last expression into the estimate for $II_{s_{m,n}}$, we get the desired estimate. So far we have finished the proof of Lemma \ref{claim1} except the Lemma \ref{lemma9}, which we will do now.\\

\noindent {\bf Proof of Lemma \ref{lemma9}:} as $x$ and $s_{m,n}$ are fixed now, later for simplicity we will just write $\tau$ instead of $\tau_{x, s_{m,n}}$. Notice that in the new coordinate we chose for $\Gamma_x$, the vector field along $\Gamma_x$ points in the direction of $(1, 0)$. Then by the definition of $\phi_s$ and $\phi_s^x$, we have
\begin{align*}
& |\phi_{s}(g_x(z),z)-\phi_{s}^x(\tau z, z)|\\
& = 2^{k-l} |\int_{\R}[\check{\psi}_0(2^{k-l}(t+g_x(z)-\tau z))-\check{\psi}_0(2^{k-l}t)]\varphi_{s}(\tau z-t, z)dt|.
\end{align*}
By the definition of the beta numbers, we have that
\begineq
|g_x(z)-\tau z|\lesim \beta_{|j_0|}(x, s_{m,n}) 2^{-k},
\endeq
which implies that
\begineq
|\check{\psi}_0(2^{k-l}(t+g_x(z)-\tau z))-\check{\psi}_0(2^{k-l}t)|\lesim 2^{-l} \beta_{|j_0|}(x, s_{m,n})
\endeq
by the fundamental theorem. In the end, non-stationary phase method leads to the final estimate:
\begin{align*}
& 2^{k-l} |\int_{\R}[\check{\psi}_0(2^{k-l}(t+g_x(z)-\tau z))-\check{\psi}_0(2^{k-l}t)]\varphi_{s}(\tau z-t, z)dt|\\
& \lesim \frac{2^{-l} \beta_{|j_0|}(x, s_{m,n}) 2^{\frac{k}{2}} 2^{\frac{k-l}{2}}}{<\min\{|m|+|n|, |m|+|n|-|j_0|\}>^N} .
\end{align*}
Thus we have finished the proof of Lemma \ref{lemma9} and hence Lemma \ref{claim1}.

Shaoming Guo, Institute of Mathematics, University of Bonn\\
\indent Address: Endenicher Allee 60, 53115, Bonn\\
\indent Email: shaoming@math.uni-bonn.de

\end{document}